%%%%%%%%%%%
% typoref.tex. V : January 18, 2000.
% Author : Anthony PHAN
% Warning : syntaxe +- LaTeX
% Sources :
% T. Lachand--Robert, ``La Ma\^\i trise de \TeX'',
% R\'ef\'erences crois\'ees;
% latex.ltx's sources;
% and of course the \TeX book.
%%%%%%%%%%%%%%%%%%%%%%%%%%%%%%%%%%%%%%%%%%%%%%%%%%%%%
%

\catcode`\@=11

%
% style (look at the behavior of \item dans \bibitem too,
% and at one ,\  in \re@dreferenceslist)
% Feel free to change:  \bibn@me (title like ``R\'ef\'erences'')
%           \bibliographym@rk (general style)
%
\def\bibn@me{R\'ef\'erences}
\def\bibliographym@rk{\centerline{{\sc\bibn@me}}
    \sectionmark\section{\ignorespaces}{\unskip\bibn@me}
    \bigbreak\bgroup
    \ifx\ninepoint\undefined\relax\else\ninepoint\fi}
%
% Beware of the \bgroup: it will be closed by \endthebibliography
%
% \refsp@ce is the spacing command that appens between multiple
% references.
%
\let\refsp@ce=\
\let\bibleftm@rk=[
\let\bibrightm@rk=]
%
% if you want more space between brackets\ldots
%\let\refsp@ce=\thinspace
%\def\bibleftm@rk{[\thinspace}
%\def\bibrightm@rk{\thinspace]}
%
% frenchy stuff
%
\def\numero{n\raise.82ex\hbox{$\fam0\scriptscriptstyle
o$}~\ignorespaces}
%
% new variables
%
\newcount\equationc@unt
\newcount\bibc@unt
\newif\ifref@changes\ref@changesfalse
\newif\ifpageref@changes\ref@changesfalse
\newif\ifbib@changes\bib@changesfalse
\newif\ifref@undefined\ref@undefinedfalse
\newif\ifpageref@undefined\ref@undefinedfalse
\newif\ifbib@undefined\bib@undefinedfalse
\newwrite\@auxout
%
% mark an equation
%
%\def\eqnum{\global\advance\equationc@unt by 1%
%\edef\lastref{\number\equationc@unt}%
%\eqno{(\lastref)}}
%
% One can reference anything, just copy the former macro
% and use it so: \machin \label{truc}
% In machin you would have defined \lastref by some number
% or any text.
%
% References macros
%
% The next macros are the core of \ref and \cite commands.
% Its first argument may be ref, pageref or bib.
%
% It is too tricky to be explained.
% (It is a bit recursive.)
% It allows using \cite or \ref or \ldots
% with arbitrary many arguments,
% for instance:
% \cite{knuth1,knuth2,ma pomme}
%
% First argument is always ref, pageref or bib.
%
\def\re@dreferences#1#2{{%
    \re@dreferenceslist{#1}#2,\undefined\@@}}
\def\re@dreferenceslist#1#2,#3\@@{\def\next{#2}%
    \expandafter\ifx\csname#1@@\meaning\next\endcsname\relax
    ??\immediate\write16
    {Warning, #1-reference "\next" on page \the\pageno\space
    is undefined.}%
    \global\csname#1@undefinedtrue\endcsname
    \else\csname#1@@\meaning\next\endcsname\fi
    \ifx#3\undefined\relax
    \else,\refsp@ce\re@dreferenceslist{#1}#3\@@\fi}
%
% notice that the former ``,\refsp@ce'' will separate
% multiple arguments. But beware of spaces
% while defining a reference or calling for it!
%
% tricky thing: \newlabel has two arguments
% {labelname}{{\lastref}{\pageref}}
% The second argument is read as two arguments
% by \newl@bel. This was necessary to get
% a jobname.aux containing the same syntax
% LaTeX would produce and use.
%
\def\newlabel#1#2{{\def\next{#1}\newl@bel#2}}
\def\newl@bel#1#2{%
    \expandafter\xdef\csname ref@@\meaning\next\endcsname{#1}%
    \expandafter\xdef\csname pageref@@\meaning\next\endcsname{#2}}
\def\label#1{{%
    \toks0={#1}\message{ref(\lastref) \the\toks0,}%
    \ignorespaces\immediate\write\@auxout%
    {\noexpand\newlabel{\the\toks0}{{\lastref}{\the\pageno}}}%
    \def\next{#1}%
    \expandafter\ifx\csname ref@@\meaning\next\endcsname\lastref%
    \else\global\ref@changestrue\fi%
    \newlabel{#1}{{\lastref}{\the\pageno}}}}
\def\ref#1{\re@dreferences{ref}{#1}}
\def\pageref#1{\re@dreferences{pageref}{#1}}
%
% bibliography macros
%
\def\bibcite#1#2{{\def\next{#1}%
    \expandafter\xdef\csname bib@@\meaning\next\endcsname{#2}}}
\def\cite#1{\bibleftm@rk\re@dreferences{bib}{#1}\bibrightm@rk}
%
% The argument of \beginthebibliography
% is any sequence of numerals which will represent
% the maximum \item's length. If you have less than 9
% \bibitem's, this argument may be {any numeral}.
% if you have between 100 and 999 \bibitem's
% this argument may be {any three numerals},
% and so on.
%
\def\beginthebibliography#1{\bibliographym@rk
    \setbox0\hbox{\bibleftm@rk#1\bibrightm@rk\enspace}
    \parindent=\wd0
    \global\bibc@unt=0
    \def\bibitem##1{\global\advance\bibc@unt by 1
        \edef\lastref{\number\bibc@unt}
        {\toks0={##1}
        \message{bib[\lastref] \the\toks0,}%
        \immediate\write\@auxout
        {\noexpand\bibcite{\the\toks0}{\lastref}}}
        \def\next{##1}%
        \expandafter\ifx
        \csname bib@@\meaning\next\endcsname\lastref
        \else\global\bib@changestrue\fi%
        \bibcite{##1}{\lastref}
        \medbreak
        \item{\hfill\bibleftm@rk\lastref\bibrightm@rk}%
        }
    }
\def\endthebibliography{\egroup\par}
%
% THE NEXT MACRO MUST BE INCLUDED
% IN THE \BYE COMMAND. FOR INSTANCE:
%
    %\catcode`@=11
    \outer\def\bye{\@closeaux
        \par\vfill\supereject\end}
    %\catcode`@=12
%
\def\@closeaux{\closeout\@auxout
    \ifref@changes\immediate\write16
    {Warning, changes in references.}\fi
    \ifpageref@changes\immediate\write16
    {Warning, changes in page references.}\fi
    \ifbib@changes\immediate\write16
    {Warning, changes in bibliography.}\fi
    \ifref@undefined\immediate\write16
    {Warning, references undefined.}\fi
    \ifpageref@undefined\immediate\write16
    {Warning, page references undefined.}\fi
    \ifbib@undefined\immediate\write16
    {Warning, citations undefined.}\fi}
%
% initialization of jobname.aux
%
\immediate\openin\@auxout=\jobname.aux
\ifeof\@auxout \immediate\write16
     {Creating file \jobname.aux}
\immediate\closein\@auxout
\immediate\openout\@auxout=\jobname.aux
\immediate\write\@auxout {\relax}%
\immediate\closeout\@auxout
\else\immediate\closein\@auxout\fi
%
% Let's read this file and open it out
%
\input\jobname.aux \par
\immediate\openout\@auxout=\jobname.aux
% this file will be closed by \bye.
%
% That's all, folks!
%

\def\bibn@me{R\'ef\'erences bibliographiques}
\catcode`@=11
\def\bibliographym@rk{\bgroup}
%
% \bye est modifie pour la biblio and la table des matieres
%
\outer\def\bye{     \par\vfill\supereject\end}

\magnification=1200

\font\tenbfit=cmbxti10
\font\sevenbfit=cmbxti10 at 7pt
\font\sixbfit=cmbxti5 at 6pt

\newfam\mathboldit

\textfont\mathboldit=\tenbfit
  \scriptfont\mathboldit=\sevenbfit
   \scriptscriptfont\mathboldit=\sixbfit

\def\bfit
{\tenbfit
   \fam\mathboldit
}

\def\Q{{\bf {Q}}}

\def\N{{\bf N}} 
\def\Z{{\bf Z}}    
\def\R{{\bf R}}

\def\Bad{{\bfit Bad}}

\def\house#1{\setbox1=\hbox{$\,#1\,$}%
\dimen1=\ht1 \advance\dimen1 by 2pt \dimen2=\dp1 \advance\dimen2 by 2pt
\setbox1=\hbox{\vrule height\dimen1 depth\dimen2\box1\vrule}%
\setbox1=\vbox{\hrule\box1}%
\advance\dimen1 by .4pt \ht1=\dimen1
\advance\dimen2 by .4pt \dp1=\dimen2 \box1\relax}

  \def\eps{{\varepsilon}}

  \def\noi{\noindent}

\def\build#1_#2^#3{\mathrel{\mathop{\kern 0pt#1}\limits_{#2}^{#3}}}

\def\date {le\ {\the\day}\ \ifcase\month\or janvier
\or fevrier\or mars\or avril\or mai\or juin\or juillet\or
ao\^ut\or septembre\or octobre\or novembre
\or d\'ecembre\fi\ {\oldstyle\the\year}}

\font\fivegoth=eufm5 \font\sevengoth=eufm7 \font\tengoth=eufm10

\newfam\gothfam \scriptscriptfont\gothfam=\fivegoth
\textfont\gothfam=\tengoth \scriptfont\gothfam=\sevengoth

\def\pro{\noindent {\it Proof.  }}

\def\smallsquare{\vbox{\hrule\hbox{\vrule height 1 ex\kern 1 ex\vrule}\hrule}}
\def\cqfd{\hfill \smallsquare\vskip 3mm}

%%%%%%%%%%%%%%%%%%%%%%%%%%%%%%%%%%%%%%%%%%%%

\def\bbQ{\bf Q}
\def\bbR{\bf R}

%%%%%%%%%%%%%%%%%%%%%%%%%%%%%%%%%%%%%%%%%%%%

\centerline{}

\vskip 4mm

\centerline{
{\bf On the complexity of a putative counterexample}}
\smallskip
\centerline{\bf to the $p$-adic Littlewood conjecture}

\vskip 8mm
\centerline{Dmitry B{\sevenrm ADZIAHIN},
Yann B{\sevenrm UGEAUD\footnote{}{\rm
2000 {\it Mathematics Subject Classification : } 11J04; 11J61, 11J83, 37A35, 37A45, 37D40.
Keywords: Diophantine approximation, Littlewood conjecture, complexity, 
continued fractions, measure rigidity.}},
Manfred E{\sevenrm INSIEDLER}\footnote{}{\rm
D.B. was supported by grant EP/L005204/1 from EPSRC.},
\& Dmitry K{\sevenrm LEINBOCK}\footnote{}{\rm
D.K. was supported in part by NSF grant DMS-1101320.}
}

{\narrower\narrower
\vskip 12mm

\proclaim Abstract. {
Let $\Vert \cdot \Vert$ denote the distance to the nearest integer
and, for a prime number $p$, let $| \cdot |_p$ denote the $p$-adic absolute
value. In 2004, de Mathan and Teuli\'e asked whether
$\inf_{q \ge 1} \, q \cdot \Vert q \alpha \Vert \cdot \vert q \vert_p = 0$
holds for every badly approximable real number $\alpha$ and every
prime number $p$. Among other results, we establish that, if the
complexity of the sequence of partial quotients of a real number
$\alpha$ grows too rapidly or too slowly, then their conjecture is true
for the pair $(\alpha, p)$ with $p$ an arbitrary prime. }

}

\vskip 11mm

\centerline{\bf 1. Introduction}

\vskip 6mm

A famous open problem in simultaneous
Diophantine approximation is
the Littlewood conjecture which claims that,
for every given pair $(\alpha, \beta)$ of real numbers, we have
$$
\inf_{q \ge 1} \, q \cdot \Vert q \alpha \Vert \cdot
\Vert q \beta \Vert = 0, \eqno (1.1)
$$
where $\Vert \cdot \Vert$ denotes the distance to the nearest integer.
The first significant contribution on this question
goes back to Cassels and Swinnerton-Dyer \cite{CaSw} who showed that
(1.1) holds when $\alpha$ and $\beta$ belong to the same cubic field.
Despite some recent remarkable progress
\cite{PoVe,EKL} the Littlewood conjecture remains an open problem.

Let ${\cal D}=(d_{k} )_{k \ge 1} $ be a sequence of integers greater
than or equal to $2$. Set $e_0 = 1$ and, for any $n \ge 1$,
$$
e_{n}  = \prod _{1 \le k \leq n}  d_{k}.
$$
For an integer $q$, set
$$
w_{{\cal D}}  (q)=\sup \{n \ge 0 : q\in  e_{n}  {\bf Z}\}
$$
and
$$
\vert q\vert _{{\cal D}}  =1/e_{w_{{\cal D}}  (q)}  =\inf
\{1/e_{n} : q\in  e_{n}  {\bf Z}\}.
$$
When $ {\cal D}$ is the constant sequence equal to $ p$, where $ p$ is a
prime number, then $ \vert \cdot \vert _{{\cal D}}$ is the usual
$ p$-adic value $| \cdot |_p$, normalized by $|p|_p = p^{-1}$.
In analogy with the Littlewood conjecture,
de Mathan and Teuli\'e \cite{BdMTe} proposed in 2004 the following conjecture.

\proclaim Mixed Littlewood Conjecture. For every real number
$\alpha$ and every sequence ${\cal D}$ of integers greater than or
equal to $2$, we have
$$
\inf_{q \ge 1} \, q \cdot \Vert q\alpha \Vert
\cdot \vert q \vert _{\cal D} =0 \eqno (1.2)
$$
holds for every real number $ \alpha  $.

Obviously, (1.2) holds if $\alpha$ is rational or
has unbounded partial quotients. Thus, we only consider the case
when $\alpha$ is an element of the set $\Bad $ of badly approximable
numbers, where
$$
\Bad  = \{ \alpha \in \R : \inf_{q \ge 1} \,
q \cdot \Vert q \alpha \Vert > 0\}.
$$
De Mathan and Teuli\'e proved that (1.2) and even the stronger
statement
$$
\liminf_{q \to + \infty} \, q \cdot  \log q \cdot \Vert q\alpha \Vert
\cdot \vert q \vert _{\cal D} < + \infty  \eqno (1.3)
$$
holds for every quadratic irrational $\alpha$ when the sequence
${\cal D}$ is bounded.

We highlight the particular case when
${\cal D}$ is the constant sequence equal to a prime number.

\proclaim $p$-adic Littlewood Conjecture. For every real number
$\alpha$ and every prime number $p$, we have
$$
\inf_{q \ge 1} \, q \cdot \Vert q \alpha \Vert \cdot
\vert q \vert_p = 0. \eqno (1.4)
$$

Einsiedler and Kleinbock \cite{EiKl07} established that, for every given prime number $p$,
the set of real numbers $\alpha$ such that the pair $(\alpha, p)$
does not satisfy (1.4) has zero Hausdorff dimension.
They also explained how to modify their proof to get an analogous result
when ${\cal D}$ is the constant sequence equal to $d \ge 2$ (not necessarily prime).

In an opposite direction, by means of a subtle
Cantor-type construction, Badziahin and Velani \cite{BaVe11}
established that, for every sequence ${\cal D}$ of integers greater
than or equal to $2$, the set of real numbers $\alpha$ such that
$$
\inf_{q \ge 3} \, q \cdot \log q \cdot \log \log q \cdot \Vert q\alpha \Vert
\cdot \vert q \vert _{\cal D}  > 0
$$
has full Hausdorff dimension. Moreover,
they showed that, for ${\cal D} = (2^{2^n})_{n \ge 1}$,
the set of real numbers $\alpha$ such that
$$
\inf_{q \ge 16} \, q \cdot \log \log q \cdot \log \log \log q \cdot \Vert q\alpha \Vert
\cdot \vert q \vert _{\cal D}  > 0
$$
has full Hausdorff dimension.

Regarding explicit examples of real numbers $\alpha$
satisfying (1.4), it was proved in \cite{BDM} that, if the sequence of partial quotients
of the real number $\alpha$ contains arbitrarily long concatenations
of a given finite block, then the pair $(\alpha, p)$ satisfies (1.4) for any prime number $p$.
A precise statement is as follows.

\proclaim Theorem BDM.
Let $\alpha  = [a_0; a_1, a_2, \ldots ]$ be in $\Bad $.
Let $T \ge 1$ be an integer and $b_1, \ldots, b_T$
be positive integers.
If there exist two sequences $(m_k)_{k \ge 1}$ and
$(h_k)_{k \ge 1}$ of positive integers with $(h_k)_{k \ge 1}$
being unbounded and
$$
a_{m_k +j + n T} = b_j, \quad
\hbox{for every $j=1, \ldots, T$ and every $n=0, \ldots, h_k-1$},
$$
then the pair $(\alpha, p)$ satisfies (1.4) for any prime number $p$.

The main purposes of the present note is to give new
combinatorial conditions ensuring that a real number satisfies
the $p$-adic (Theorem 2.1) and the mixed (Corollary 2.4) Littlewood conjectures
and to study the complexity of the continued fraction expansion of a
putative counterexample to (1.2) or (1.4); see Theorem 2.1 and
Corollary 2.4 below. Furthermore, in Section 3 we make a connection
between the mixed Littlewood conjecture and a problem on the
evolution of the sequence of the Lagrange constants of the multiples
of a given real number. Proofs of our results are given in Sections
4 to 6.

Throughout the paper,
we assume that the reader is familiar with the classical results from the theory
of continued fractions.

\vskip 5mm

\goodbreak

\centerline{\bf 2. New results on the mixed and
the $p$-adic Littlewood conjectures}

\vskip 6mm

To present our results, we adopt a point of view from
combinatorics on words. We look at the continued fraction
expansion of a given real number $\alpha$
as an infinite word.

For an infinite word ${\bf w} = w_1 w_2 \ldots$
written on a finite alphabet
and for an integer $n \ge 1$, we denote
by $p(n, {\bf w})$ the number of distinct blocks of $n$ consecutive
letters occurring in ${\bf w}$, that is,
$$
p(n, {\bf w}) := {\rm Card} \{w_{\ell+1} \ldots w_{\ell + n} : \ell
\ge 0\}.
$$
The function $n \mapsto p(n, {\bf w})$ is called the {\it complexity
function} of ${\bf w}$. For
a badly approximable real number $\alpha =
[a_0; a_1, a_2, \ldots]$, we set
$$
p(n, \alpha) := p(n, a_1 a_2 \ldots), \quad n \ge 1,
$$
and we call $n \mapsto p(n, \alpha)$ the {\it complexity
function} of $\alpha$.
Observe that, for all positive integers $n, n'$, we have
$$
p(n + n', \alpha) \le p(n, \alpha) \cdot p(n', \alpha),
$$
thus, the sequence $(\log p(n, \alpha))_{ n \ge 1}$ is subadditive
and  the sequence $((\log p(n, \alpha))/n)_{ n \ge 1}$ converges.

In the present paper we show that if the real number $\alpha$ is a
counterexample to the $p$-adic Littlewood conjecture, then its
complexity function $n \mapsto p(n, \alpha)$ can neither increase
too slowly nor too rapidly as $n$ tends to infinity.

\bigskip
\goodbreak

\centerline{\bf 2.1. High complexity case}

\vskip 4mm

For a positive integer $K$, set
$$
\Bad_K:= \{\alpha = [a_0; a_1, a_2\ldots]\;:\; a_i\le K, i \ge 1 \}
$$
and observe that the set of badly approximable numbers is the union over
all positive integers $K$ of the sets $\Bad_K$.
It immediately follows from the definition of the
complexity function $n \mapsto p(n, \alpha)$ that, for
every $\alpha$ in $\Bad_K$ and every $n \ge 1$, we have
$$
p(n, \alpha) \le K^n.
$$
Consequently, the complexity function of the continued fraction of any number
$\alpha$ in $\Bad$ grows at most exponentially fast. Our first result shows that a
putative counterexample to the $p$-adic Littlewood conjecture must
satisfy a much more restrictive
condition.

\proclaim Theorem 2.1.
Let $\alpha$ be a real number satisfying
$$
\lim_{n\to\infty}{\log p(n,\alpha) \over n} >0.  \eqno (2.1)
$$
Then, for every prime number $p$, we have
$$
\inf_{q \ge 1} q\cdot ||q\alpha||\cdot |q|_p=0.
$$

In other words the complexity of the continued fraction expansion of every potential
counterexample to the $p$-adic Littlewood conjecture must grow subexponentially.

Our proof relies on a~$p$-adic generalisation of the measure classification
result in~\cite{Lindenstrauss} (provided by~\cite{EinLin}), the connection
between such dynamical results and the Diophantine approximation problem
as was used before in~\cite{EKL,EiKl07}, and the observation that one counterexample
actually gives rise to many more counterexamples (see Proposition 4.1).

\bigskip

\centerline{\bf 2.2. Low complexity case}

\vskip 4mm

A well-known result of Morse and Hedlund \cite{MoHe38,MoHe40}
asserts that $p(n, {\bf w}) \ge n + 1$ for $n \ge 1$, unless ${\bf
w}$ is ultimately periodic (in which case there exists a constant
$C$ such that $p(n, {\bf w}) \le C$ for $n \ge 1$). Infinite words
${\bf w}$ satisfying $p(n, {\bf w}) = n + 1$ for every $n \ge 1$ do
exist and are called {\it Sturmian words}. In the present paper we show
that if $\alpha$ is a counterexample to the $p$-adic (or, even, to the mixed)
Littlewood conjecture, then the lower bound for the
complexity function of $\alpha$ must be stronger than this estimate.
Before stating our result we give a classical definition (see e.g. \cite{AlSh}).

\proclaim Definition 2.2.
An infinite word ${\bf w}$ is
recurrent if every finite block occurring in ${\bf w}$ occurs
infinitely often.

Classical examples of recurrent infinite words include periodic words, Sturmian
words, the Thue--Morse word, etc.

\proclaim Theorem 2.3.
Let $(a_k)_{k \ge 1}$ be a sequence of
positive integers. If there exists an integer $m \ge 0$
such that the infinite word $a_{m+1} a_{m+2} \ldots $ is recurrent,
then, for every sequence ${\cal D}$ of integers greater
than or equal to $2$, the real number $\alpha := [0; a_1, a_2, \ldots ]$
satisfies
$$
\inf_{q \ge 1} \, q \cdot \Vert q\alpha \Vert
\cdot \vert q \vert _{\cal D} =0.
$$

As a particular case, Theorem 2.3 asserts that (1.2) holds for every
quadratic number $\alpha$ and every sequence ${\cal D}$ of integers
greater than or equal to $2$, including unbounded sequences (unlike
in \cite{BdMTe}, where ${\cal D}$ is assumed to be bounded). Unlike
in \cite{BdMTe}, our proof does not use $p$-adic analysis.

Theorem 2.3 implies a non-trivial lower bound for the complexity
function of the continued fraction expansion of a putative
counterexample to (1.2).

\proclaim Corollary 2.4. Let ${\cal D}$ be a sequence of integers
greater than or equal to $2$ and $\alpha$ be a real number such that
the pair $(\alpha, {\cal D})$ is a counterexample to the mixed
Littlewood conjecture (i.e., does not satisfy (1.2)). Then, the
complexity function of $\alpha$ satisfies
$$
\lim_{n \to + \infty} \, p(n, \alpha) - n = + \infty.
$$

The next corollary highlights a special family of infinite recurrent words.
A finite word $w_1 \ldots w_n$ is called a {\it palindrome}
if $w_{n+1-i} = w_i$ for $i=1, \ldots , n$.

\proclaim Corollary 2.5.
Let $(a_k)_{k \ge 1}$ be a sequence of
positive integers. If there exists an increasing sequence $(n_j)_{j
\ge 1}$ of positive integers such that $a_1 \ldots a_{n_j}$ is a
palindrome for $j \ge 1$, then, for every sequence ${\cal D}$ of
integers greater than or equal to $2$, the real number $\alpha :=
[0; a_1, a_2, \ldots ]$ satisfies
$$
\inf_{q \ge 1} \, q \cdot \Vert q\alpha \Vert
\cdot \vert q \vert _{\cal D} =0.
$$

As shown in Section 6, our approach allows us to give an alternative
proof to (1.3) when $\alpha$ is quadratic irrational and ${\cal D}$
is bounded. Furthermore, we are able to quantify Theorem 2.3 for a
special class of recurrent words.

\proclaim Definition 2.6.
We say that an infinite word ${\bf w}$ is linearly
recurrent if there exists $C > 1$ such that the
distance between two consecutive occurrences of any finite block $W$
occurring in ${\bf w}$ is bounded by $C$ times the length of $W$.

We obtain the following quantitative result.

\proclaim Theorem 2.7.
Let $(a_k)_{k \ge 1}$ be a bounded sequence of positive integers.
If there exists an integer $m \ge 0$
such that the infinite word $a_{m+1} a_{m+2} \ldots $ is linearly recurrent,
then, for every sequence ${\cal D}$ of integers greater
than or equal to $2$, the real number $\alpha := [0; a_1, a_2, \ldots ]$
satisfies
$$
\liminf_{q \to + \infty} \, q \cdot (\log \log q)^{1/2} \cdot \Vert q\alpha \Vert
\cdot \vert q \vert _{\cal D} < + \infty.
$$

\bigskip

\goodbreak

\centerline{\bf 2.3. Comparison with the Littlewood conjecture}

\vskip 4mm

According to Section 5 of \cite{BdMTe}, the initial motivation of the introduction of the mixed
Littlewood conjecture was the study of a problem quite close to the 
Littlewood conjecture, but seemingly a little simpler,
with the hope to find new ideas suggesting a possible approach towards the
resolution of the Littlewood conjecture itself.

We are not aware of any relationship between both conjectures. For instance, a real
number $\alpha$ being given, we do 
not know any connection between the two statements 
`(1.2) holds for every sequence ${\cal D}$' and 
`(1.1) holds for every real number $\beta$'.

The interested reader is directed to \cite{Bu14} for a survey of 
recent results and developments on and around the Littlewood conjecture
and its mixed analogue. He will notice that the state-of-the-art regarding the
Littlewood and the $p$-adic Littlewood conjectures is essentially, but
not exactly, the same.

For instance, Theorem 5 in \cite{Lin10} asserts that for every real number $\alpha$
with (2.1), we have
$$
\inf_{q \ge 1} \, q \cdot \| q \alpha \| \cdot
\| q \beta \| = 0,
$$
for every real number $\beta$. This is the exact analogue to Theorem 2.1 above. 
However, the low complexity case remains very mysterious for the
Littlewood conjecture, since we even do not know whether or not it holds
for the pair $(\sqrt{2}, \sqrt{3})$.

\vskip 5mm

\goodbreak

\centerline{\bf 3. On the Lagrange constants of the multiples of a real number}

\vskip 6mm

Our main motivation was the study of the $p$-adic and the mixed Littlewood
conjectures. However, the proofs of Theorems 2.3 and 2.8
actually give us much stronger results on the behaviour of the Lagrange
constants of the multiples of certain real numbers.

\proclaim Definition 3.1.
The Lagrange constant $c(\alpha)$
of an irrational real number $\alpha$ is the quantity
$$
c(\alpha) := \liminf_{q \to + \infty} \, q \cdot ||q \alpha||.
$$

Clearly, $\alpha$ is in $\Bad$ if and only if $c(\alpha) > 0$. A
classical theorem of Hurwitz (see \cite{Per,BuLiv}) asserts that
$c(\alpha) \le 1/ \sqrt{5}$ for every irrational real number
$\alpha$.

For any positive
integer $n$ and any badly approximable number $\alpha$ we have
$$
{c(\alpha) \over n} \le c(n \alpha ) \le n c(\alpha).   \eqno (3.1)
$$
To see this, note that
$$
\Bigl| n \alpha - {np \over q} \Bigr| = n \Bigl| \alpha - {p \over q} \Bigr|
$$
and
$$
\Bigl|  \alpha - {p \over n q} \Bigr| = {1 \over n} \,  \Bigl| n \alpha - {n p \over n q} \Bigr|.
$$

The first general result on the behaviour of the sequence $(c(n \alpha))_{n \ge 1}$
is Theorem 1.11 of Einsiedler, Fishman, and Shapira \cite{EFS}, reproduced below.

\proclaim Theorem EFS.
Every badly approximable real number $\alpha$ satisfies
$$
\inf_{n \ge 1} \, c(n \alpha) = 0.
$$

Theorem EFS motivates the following question.

\proclaim Problem 3.2.
Prove or disprove that every badly approximable real number $\alpha$ satisfy
$$
\lim_{n \to + \infty} \, c(n \alpha) = 0.  \eqno (3.2)
$$

There is a clear connection between Problem 3.2 and the mixed
Littlewood conjecture. Indeed, if $\alpha$ satisfies (3.2) and
if ${\cal D}$ is as in Section 1, then, keeping the notation from this
section, for every $\eps > 0$, there exists a positive integer $n$
such that $c(e_n \alpha) < \eps$. Consequently, there are
arbitrarily large integers $q$ with the property that
$$
q \cdot || q e_n \alpha || < \eps
$$
thus,
$$
q e_n \cdot || q e_n \alpha || \cdot |q e_n|_{\cal D} < \eps ,
$$
since $|q e_n|_{\cal D} \le 1/e_n$. This proves that (1.2)
holds for the pair $(\alpha, {\cal D})$.

Our proof of Theorem 2.3 actually gives the following
stronger result.

\proclaim Theorem 3.3.
Let $(a_k)_{k \ge 1}$ be a sequence of positive integers.
If there exists an integer $m \ge 0$
such that the infinite word $a_{m+1} a_{m+2} \ldots $ is recurrent,
then the real number $\alpha := [0; a_1, a_2, \ldots ]$ satisfies (3.2) and, moreover,
$$
c(n \alpha) \le {8 q_m^2 \over n}, \quad \hbox{for $n \ge 1$,}
$$
where $q_m$ denotes the denominator of the rational
number $[0; a_1, \ldots , a_m]$.

In view of the left-hand inequality of (3.1), the conclusion
of Theorem 3.3 is nearly best possible.  

Using the same arguments as for the proof of Corollary 2.4, we
establish that the complexity function of a real number which does
not satisfy (3.2) cannot be too small.

\proclaim Corollary 3.4.   
Let $\alpha$ be a real number such that
$$
\sup_{n \ge 1} \, n \, c(n \alpha) = + \infty.
$$
Then, the complexity function of $\alpha$ satisfies
$$
\lim_{n \to + \infty} \, p(n, \alpha) - n = + \infty.
$$

\vskip 5mm

\goodbreak

\centerline{\bf 4. High complexity case}

\vskip 5mm

We follow the interpretation of the $p$-adic Littlewood conjecture used
by Einsiedler and Kleinbock in~\cite{EiKl07} and consider the
following more general problem:

\proclaim Generalized $p$-adic Littlewood Conjecture.
For every prime
number $p$ and for every pair $(u,v)\in {\bbR_{>0}}\times {{\bbQ}_p}$ we have
$$
\inf_{a\in\N,b\in \N\cup\{0\}} \max\{|a|,|b|\}\cdot |au-b|\cdot |av-b|_p
=0. \eqno(4.1)
$$

Clearly~$\alpha$ satisfies the $p$-adic Littlewood conjecture
(i.e.~(1.4)) if and only if~$-\alpha$ satisfies the $p$-adic
Littlewood conjecture. For that reason we restrict our attention to
positive numbers. Moreover, one can check (see for
example~\cite{EiKl07}, a discussion after Theorem 1.2) that if
$\alpha$ is a counterexample to the $p$-adic Littlewood conjecture
then $(\alpha^{-1},0)$ is a counterexample to the above generalized
$p$-adic Littlewood conjecture. The next proposition goes further
and shows that one counterexample $\alpha$ to the $p$-adic
Littlewood conjecture provides a countable collection of
counterexamples to the generalized $p$-adic Littlewood conjecture.

\proclaim Proposition 4.1. Let $p$ be a prime number and $\alpha>0$
an irrational number. Let $\eps$ be in $(0, 1/2]$ and assume that
$$
\inf_{q\ge 1} q\cdot ||q\alpha||\cdot |q|_p>\eps.
$$
Then, we have
$$
\inf_{a\in\N,b\in\N\cup\{0\}}\;\; \max\{a,b\} \cdot \left|a\cdot
{||q_n\alpha||\over ||q_{n-1}\alpha||}-b\right|\cdot \left|a\cdot
\left({q_n \over q_{n-1}}\right)+b\right|_p>{\eps^2 \over 4},
\eqno(4.2)
$$
where $(q_k)_{k \ge 1}$ is the sequence of the
denominators of the convergents to $\alpha$.

Note that, writing $\alpha = [a_0; a_1, a_2, \ldots]$, we have
$$
{q_n \over q_{n-1}} = [0;a_n,a_{n-1},\ldots,a_1]\quad{\rm and}\quad
{||q_n\alpha||\over ||q_{n-1}\alpha||} = [0;a_{n+1},a_{n+2},\ldots]\in(0,1),
$$
for every $n \ge 1$.

\noi {\it Proof of Proposition~4.1.} We assume that
$||q\alpha||\cdot |q|_p>\eps / q$ for every integer $q \ge 1$. We
use the classical estimate from the theory of continued fractions
$$
||q_n\alpha||<q_{n+1}^{-1}<(a_nq_n)^{-1}.
$$
It implies that $q_n||q_n\alpha||< a_n^{-1}$ and hence
$a_n<\eps^{-1}$. In other words
$$
\alpha\in \Bad_N,\quad{\rm where}\; N:= [\eps^{-1}]\quad {\rm
and}\quad q_n<(N+1)q_{n-1}.
$$

Now choose some~$a\geq 1$,~$b\geq 0$, and modify the left hand side of~(4.2):
$$
\max\{a,b\} \cdot \left|a {||q_n\alpha||\over
||q_{n-1}\alpha||}-b\right|\cdot \left|a \left({q_n \over
q_{n-1}}\right)+b\right|_p
$$$$
={\max\{a,b\}\over ||q_{n-1}\alpha||\cdot |q_{n-1}|_p}\cdot \big|a
||q_n\alpha||-b||q_{n-1}\alpha||\;\big|\cdot |aq_n+bq_{n-1}|_p.
$$
Since $||q_{n-1}\alpha||<(a_{n-1}q_{n-1})^{-1}\le q_{n-1}^{-1}$, the
first term in this product is bounded from below by
$\max\{a,b\}\cdot q_{n-1}$, which in turn is $\geq (N+1)^{-1}
\cdot\max\{aq_n, bq_{n-1}\}$. The second term is estimated as
follows:
$$
\big|a\cdot ||q_n\alpha||-b\cdot||q_{n-1}\alpha||\;\big| =
|(aq_n+bq_{n-1})\alpha -(ap_n+bp_{n-1})|\ge
||(aq_n+bq_{n-1})\alpha||.
$$
Therefore, a lower bound of the whole product is
$$
{\max\{aq_n, bq_{n-1}\}\over (N+1)} \cdot
||(aq_n+bq_{n-1})\alpha||\cdot |aq_n+bq_{n-1}|_p
$$$$
\ge{\eps\over N+1}\cdot
{\max\{a q_n,b q_{n-1}\}\over aq_n+bq_{n-1}}\ge {\eps\over 2(N+1)}.
$$
Since $N+1 \le \eps^{-1}+2$ and $\eps \le 1/2$, this proves the proposition. \cqfd

In~\cite{EiKl07} the authors showed that the set of counterexamples to the
generalized $p$-adic Littlewood conjecture is rather small.
More precisely, the following theorem
(Theorem 5.2 in \cite{EiKl07}) was stated there, along with a scheme of proof.

\proclaim Theorem EK.
Let $p$ be a prime number.
Then the set of pairs $(u,v)\in {\bbR}\times {{\bbQ}_p}$ which do not satisfy
$$
\liminf_{a,b\in \Z} |a|\cdot |au-b|\cdot |av-b|_p
=0.\eqno(4.3)
$$
is a countable union of sets of box dimension zero.

The outlined proof was based on a theorem
due to Einsiedler and Lindenstrauss \cite{EinLin}
which at the time of publication of \cite{EiKl07} had not appeared yet.
In the present paper  we reexamine the methods of \cite{EiKl07} and
provide a more precise result regarding the set of pairs $(u,v)\in [0,1]\times \Z_p$ which do not
 satisfy~(4.1).

\proclaim Theorem 4.2.
For every prime number $p$,   
the set of pairs $(u,v)\in [0,1]\times \Z_p$ which do not satisfy (4.1)
is a countable union of sets of box dimension zero. Moreover,
for every $\eps>0$ the set of $(u,v)\in [0,1]\times \Z_p$ which satisfy
$$
\inf_{a\in \N,b\in\N\cup\{0\}} \max\{a,b\}\cdot |au-b|\cdot |av-b|_p \ge
\eps.\eqno(4.4)
$$
has box dimension zero.

The proof of Theorem 4.2 is postponed to Section 5.

Roughly speaking, this theorem together with Proposition~4.1 implies that for every
counterexample $\alpha = [0; a_1, a_2, \ldots ]$ to the $p$-adic Littlewood conjecture the set
$$
\{[0;a_m,a_{m+1},\ldots]\;:\;m \ge 1 \}
$$
has box dimension zero.   
Let us now show that this in turn  implies the
statement of Theorem~2.1.

\vskip 4mm

\noi {\it Proof of Theorem~2.1.}
We will prove the contrapositive of the theorem.
So assume that~$\alpha$ is a counterexample to (1.4).
By the homogeneity of  (1.4)
 we may assume that~$\alpha=[a_0;a_1,\ldots]$ is positive.
By Proposition~4.1 this leads
to a countable collection
$$
 B=\Bigl\{ (\alpha_n,\beta_n)
 = \Bigl( [0;a_{n+1},a_{n+2},\ldots],{q_n \over q_{n-1}}\Bigr)
 : n\geq 1\Bigr\}
$$
of pairs in~$[0,1]\times\Q_p$ that all satisfy (4.2).  By Theorem~4.2,
the set
$$
B\cap [0,1]\times\Z_p=\{(\alpha_n,\beta_n):n\geq 1 \enspace
\hbox{and~$|q_{n-1}|_p=1$}\}
$$
has box dimension zero.

Let~$N\geq 1$ be such that~$a_n\in\{1,\ldots,N\}$ for all~$n\geq 1$.
We set~
$$
S'=\{\ell\geq 1: |q_{\ell-1}|_p=1\}.
$$
Let~$\delta>0$ be arbitrary and let~$\pi_\infty:[0,1]\times\Z_p\to[0,1]$ denote the projection
to the real coordinate.
Then, the definition of box dimension shows that, for all sufficiently large~$n$, the set
$$
B'=\{\alpha_\ell:\ell\in S'\}=\pi_\infty(B\cap [0,1]\times\Z_p)
$$
can be covered by~$s_n\leq e^{\delta n}$ intervals~$I_1,\ldots,I_{s_n}$
of size~$(1+N)^{-2n}$.

We also define another disjoint collection of intervals.
To any word~$w = w_1 \ldots w_n$ in $\{1,\ldots,N\}^n$,
we associate the interval $[w]$ composed of the real
numbers in $(0, 1)$ whose first $n$ partial quotients are $w_1, \ldots , w_n$. 
The basic
properties of continued fractions show that the length of~$[w]$ is at most~$2^{-n+2}$
and at least~$(1+N)^{-2n}$. It follows that a given interval~$I_j$ from the above list
can intersect at most two (neighbouring) intervals of the form~$[w]$ for~$w\in\{1,\ldots,N\}^n$.
This implies that
$$
{\rm Card} \{a_{\ell+1}\ldots a_{\ell+n}:\ell \in S'\} \leq 2 s_n\leq 2 e^{\delta n}.
$$
To remove the restriction~$\ell\in S'$ in the above counting, we note that~$\ell\notin S'$
implies~$\ell+1\in S'$ since~$q_{\ell-1}$ and~$q_\ell$ are coprime, by the properties
of continued fractions. Therefore,
$$
p(n,\alpha)= {\rm Card} \{a_{\ell+1}\ldots a_{\ell+n}:\ell \geq 0\}
\leq 2 s_n+2 N s_{n-1}\leq 2(1+N) e^{\delta n}.
$$
As~$\delta>0$ was arbitrary, the theorem follows.
  \cqfd

\vskip 4mm

\centerline{\bf 5. Measure Rigidity and the Proof of Theorem 4.2.}

\vskip 4mm

We follow the strategy
outlined in \cite{EiKl07} (which in turn  generalizes the
argument  from \cite{EKL}). For this,
we set
$$
G={\rm SL}_2(\R)\times {\rm SL}_2(\Q_p), \ \Gamma={\rm SL}_2(\Z[1/p]),
\ X =G/\Gamma\,,\eqno(5.1)
$$
where $\Z[1/p]$ is embedded diagonally via~$a\mapsto (a,a)$ in $\R\times\Q_p$.
 In other words, for  $(A,B) \in  G$,
points $x=(A,B)\Gamma\in X$ are identified with
unimodular 
lattices $(A,B)\Z[1/p]$ in $\R^2\times\Q_p^2$ that are
generated by the column vectors of $A$ and $B$.

We also set
$$
\psi(t,n)=\left (\left(\matrix{ e^{-t}&0\cr 0&e^t }\right),
\left(\matrix{p^n&0\cr 0&p^{-n }}\right)\right)\eqno(5.2)
$$
for~$(t,n)\in\R\times\Z$, and define the cone
$$
C=\{(t,n)\mid n\geq 0, e^t p^{-n}\geq 1\}\,.\eqno(5.3)
$$
Furthermore, for~$(u,v)\in\R\times\Q_p$,
we define the coset (which we will think of as a point)
$$
x_{u,v}=\left (\left(\matrix{ 1&0\cr u&1 }\right),
\left(\matrix{1&0\cr v&1}\right)\right)\Gamma.
$$

Compact subsets of $X$ can be characterized by the analogue
of Mahler's compactness criterion
(see \cite{EiKl07}, Theorem\ 2.1) so that a subset $K\subset X$
has compact closure if and only if
there exists some~$\delta>0$  so that $K\subset K_\delta$, with
$$
 K_\delta=\left\{ g\Gamma\in K  :  g\Z[1/p]^2\cap B_\delta^{\R^2\times\Q_p^2}=\{0\}\right\},
$$
where $B_\delta$ denotes the ball of radius $\delta$ centered at
zero.

In \cite{EiKl07}, Proposition\ 2.2, a connection between unboundedness of the
cone orbit
$$
\psi(C)x_{u,v}=\{\psi(t,n)x_{u,v}  : (t,n)\in C\}
$$
and (4.3) is given. However, we will need to show the following
refinement.

\proclaim Proposition 5.1.
Let $(u,v)\in 
(0,1)\times\Z_p$ and $0 < \varepsilon<1$ be arbitrary. If $(u,v)$ satisfies (4.4),
then  
$\psi(C) x_{u,v}\subset K_\delta$ for~$\delta =(\eps/2)^{1/3}$.
\vskip 4mm

We note that in \cite{EiKl07} the converse of the above implication,
in a slightly different form, has also been claimed (without a proof).
But the other direction is not clear and luckily is also not needed for the proof of our results (or
the results of~\cite{EiKl07}).
\vskip 4mm

\noi {\it Proof of Proposition~5.1.}
Take $\delta =(\eps/2)^{1/3}$ and suppose that $\psi(C) x_{u,v}$
is not contained in  $K_\delta$; that is, there exists a pair $(t,n)$ with $n\ge 0$
  and $e^tp^{-n}\geq 1$ such that $\psi(t,n)x_{u,v}\Z[1/p]$
  contains a nonzero element in $B_\delta^{\R^2\times\Q_p^2}$.

  Clearly, $\psi{(t,n)}x_{u,v}$ is generated by
$$
   \left(\pmatrix {e^{-t}\cr  e^tu},
   \pmatrix {p^{n}\cr  p^{-n}v}\right){ \rm and }
   \left(\pmatrix{ 0\cr  e^t},
   \pmatrix{ 0\cr p^{-n}}\right).
$$
  However, since $\psi{(t,n)}x_{u,v}$ is a
  $\Z[1/{p}]$-module, the vectors
$$
   \left(\pmatrix{ e^{-t}p^{-n}\cr  e^tp^{-n}u},
   \pmatrix{ 1\cr  p^{-2n}v}\right){ \rm and }
   \left(\pmatrix{ 0\cr  e^tp^{-n}},
   \pmatrix{ 0\cr p^{-2n}}\right)
$$
  are also generators. Therefore, there exists some nonzero
  $(a,b)\in\Z[1/{p}]^2$ such that
$$
   \left(\pmatrix{e^{-t}p^{-n}a\cr  e^tp^{-n}(au-b)},
   \pmatrix{a\cr p^{-2n}(av-b)}\right)\in\R^2\times\Q_p^2
$$
  is $\delta$-small. In particular,
$
|a|_p$ is less than $\delta$,
which  
implies that $a\in\Z$. Since
  $n\geq 0$ and $v\in\Z_p$, the inequality
$$
p^{2n}|av-b|_p = |p^{-2n}(av-b)|_p <\delta    \eqno (5.4)
$$
shows that  $b\in\Z$ as well.  
Also, since $|u| < 1$, the inequalities $t \ge 0$,
$|e^{-t}p^{-n}a| <\delta$
and
$$
|e^{t}p^{-n}(au-b)| <\delta   \eqno (5.5)
$$
imply that
$$
e^{-t}p^{-n}\max\{ |a|,|b| \} <2\delta.  \eqno (5.6)
$$
By taking the product of the inequalities  (5.4), (5.5) and (5.6),
we arrive at
$$
\max \{ |a|,|b| \} \cdot|au-b|\cdot|av-b|_p<2\delta^3 = \eps.
$$
Also note that~$u>0$,~(5.5) and~$e^tp^{-n}\geq 1$ imply that~$a$ and~$b$
  have the same sign (in the sense that~$ab\geq 0$).
  Without loss of generality we may assume~$a,b\geq 0$.
  Similarly,~$b=0$ implies~$a=0$ and contradicts our choice of~$(a,b)$.
  However,~$a\geq 1$ and~$b\geq 0$ contradicts (4.4).
Consequently, $\psi(C) x_{u,v}$
is contained in  $K_\delta$.
\cqfd
\vskip 4mm

We also need the following partial measure classification result.

\proclaim Theorem 5.2.
 The Haar measure is the only~$\psi$-invariant and ergodic probability measure $\mu$
  on~$X$ for which some~$(t,n)\in\R\times\Z$ has positive entropy $h_\mu(\psi(t,n))>0$.
 \vskip 4mm

Theorem 5.2 follows from Theorem 1.3 of \cite{EinLin}, recalled below.
We write~$\infty$ for the archimedean place of~$\Q$ and~$\Q_\infty=\R$. Moreover,
for a finite set~$S$ of places we define~$\Q_S=\prod_{v\in S}\Q_v$
for the corresponding product of the local fields.

\proclaim Theorem EL.
    Let $\bf G$ be a $\Q$-almost simple linear algebraic group, let~$S$
    be a finite set of places containing the archimedean place~$\infty$,
    let $\Gamma<G={\bf G} (\Q_S)$ be an arithmetic lattice, and let~$X=G/\Gamma$.
    Finally, let~$A$ be the direct
    product of maximal~$\Q_v$-diagonalizable algebraic subgroups
    of~${\bf G}(\Q_v)$ for~$v\in S$.
    Let $\mu$ be an $A$-invariant and ergodic probability measure on $X$.
    Suppose in addition that
    $\mu$ is not supported on any periodic orbit
    $ g{\bf L}(\Q_S)\Gamma $ for any 
$g\in G$ and proper reductive $\Q$-subgroup ${\bf  L}<\bf G$,
      that ${\rm rank}(A)\geq 2$, and that $h_\mu(a)>0$  for some $a\in A$.
    Then there is a finite index subgroup $L<G$ so that $\mu$ is
    $L$-invariant and supported on a single $L$-orbit.

\noindent {\it Proof of Theorem 5.2.} We let $G, \Gamma, X$ be as in (5.1).
In the special case~$(t,n)=(t,0)$ this follows directly from \cite{Lindenstrauss}, Theorem~1.1.
The method of proof of~\cite{Lindenstrauss}, Theorem~1.1.~would in principle also
give the general case of Theorem~5.2, but we  will instead derive it from the more
general Theorem EL.

Assume now~$\mu$ is a~$\psi$-invariant and ergodic probability
measure with positive entropy for some~$(t,n)\in\R\times\Z$.
Strictly speaking~${\rm Im}(\psi)=\psi(\R\times\Z)$ does not equal
the product~$A$ of the full diagonal subgroup of~${\rm SL}(\R)$ and
the full diagonal subgroup of~${\rm SL}(\Q_p)$. However,~$K=A/{\rm
Im}(\psi)$ is compact which allows us to define the~$A$-invariant
and ergodic measure~$\mu_A=\int_{K}a_*\mu {\rm d} (a{\rm Im}(\psi))$
with positive entropy for~$\psi(t,n)$.

Note that a proper nontrivial reductive subgroup~$\bf L$ of~${\rm SL}_2$
must be a diagonalisable subgroup. However,
if~$\mu_A$ would be supported on a single orbit 
$g{\bf L}(\Q_S)\Gamma$ this would mean that~$\mu_A$
is supported on a single periodic orbit for~$A$ and would force entropy to be equal to zero.
Therefore, all assumptions to Theorem EL are satisfied and it follows that~$\mu_A$ is
invariant under a finite index subgroup of
$G$. However,
$G$ does not have any proper
finite index subgroup and so~$\mu_A$ must be the Haar measure on~$X$.
The definition of~$\mu_A$
now expresses the Haar measure as a convex combination of~$\psi$-invariant measures.
By ergodicity
of the Haar measure under the action of~$\psi(\R\times\Z)$ this implies that~$\mu$ 
equals the Haar measure on $X$ also.
\cqfd

Finally, for the proof of Theorem~4.2, we need to quote another result
highlighting a connection between entropy and box dimension.
Recall that given $g\in G$, the {\sl unstable
 horospherical subgroup for~$g$\/} is  the maximal subgroup of $G$
 such that each of its elements $h$
 satisfies~$g^{-j}h g^j\to 1$ as $j\to\infty$.  
 What follows is a special case of  Proposition 4.1 from \cite{EiKl07},
 cf.\ also Proposition 9.1 from \cite{EKL}:

\proclaim Proposition EK.
    Let   $G, \Gamma, X$  be as in (5.1),  $\psi$ be as in (5.2) and $C$ as in (5.3).
    Take $(t,n)\in C$ and let $Y\subset X$ be a compact set such that
    no $\psi$-invariant and ergodic
  probability measure supported on $Y$ has positive entropy for
$\psi(t,n)$.
  Then for any compact subset $B$ of the unstable horospherical subgroup for $\psi(t,n)$
  and any $x\in X$  the set
 $$
    \left\{u\in B:\psi(C) ux\subset Y\right\}
$$
   has box dimension zero.

\noi {\it Proof of Theorem~4.2.}
  Fix some~$\varepsilon>0$,  
let $\delta=\varepsilon^{1/3}$ and $Y=K_\delta\subset X$.
Also pick~$t > 0$ such that $(t,1)\in C$.
 Note that~$\{x_{u,v} : (u,v)\in\R\times\Q_p\}$
 is the unstable
 horospherical subgroup
 for~$\psi(t,1)$.
  By Theorem 5.2
  and since~$Y$ is a proper closed subset of~$X$,
  there is no~$\psi$-invariant and ergodic probability measure supported on~$Y$,
  which is precisely the assumption  
  of Proposition EK.
 By that result we obtain that the set of~$(u,v)\in [0,1]\times\Z_p$
 with~$\psi(C)x_{u,v}\subset Y$ has box dimension zero.
 However, by Proposition 5.1, this implies
 that
 the set of $(u,v)\in (0,1)\times\Z_p$  satisfiying (4.4) has box dimension zero, 
 and  Theorem~4.2 follows.
\cqfd
\vskip 4mm

\vskip 4mm

\centerline{\bf 6. Low complexity case}

\bigskip

\centerline{\bf 6.1. Auxiliary results}

\vskip 6mm

We begin with two classical lemmata on continued fractions, whose
proofs can be found for example in Perron's book \cite{Per}.

For positive integers $a_1, \ldots, a_n$, denote
by $K_n (a_1, \ldots, a_n)$ the denominator of the rational number
$[0; a_1, \ldots, a_n]$. It is commonly called a {\it continuant}.

\proclaim Lemma 6.1. For any positive integers $a_1, \ldots, a_n$
and any integer $k$ with $1 \le k \le n-1$, we have
$$
K_n (a_1, \ldots, a_n)  = K_n (a_n, \ldots, a_1)
$$
and
$$
\eqalign{
K_k (a_1, \ldots, a_k) \cdot K_{n-k} (a_{k+1}, \ldots, a_n)
& \le K_n (a_1, \ldots , a_n) \cr
& \le 2 \, K_k (a_1, \ldots, a_k) \cdot K_{n-k} (a_{k+1}, \ldots, a_n). \cr}
$$

\proclaim Lemma 6.2. Let $\alpha = [0; a_1, a_2, \ldots]$ and $\beta
= [0; b_1, b_2, \ldots]$ be real numbers. Assume that there exists a
positive integer $n$ such that $a_i = b_i$ for any $i=1, \ldots, n$.
We then have $|\alpha - \beta| \le K_n (a_1, \ldots , a_n)^{-2}$.

A homogeneous linear recurrence sequence with constant
coefficients ({\it recurrence sequence} for short) is a
sequence $(u_n)_{n \ge 0}$ of complex numbers such that
$$
u_{n+d} = v_{d-1} u_{n+d-1} + v_{d-2} u_{n+d-2} +
\ldots + v_0 u_n \quad (n \ge 0),
$$
for some complex numbers $v_0, v_1, \ldots , v_{d-1}$
with $v_0 \not=0$ and with initial values $u_0, \ldots , u_{d-1}$
not all zero. The positive integer $d$ is called the
{\it order} of the recurrence.

\proclaim Lemma 6.3. Let $(u_n)_{n \ge 1}$ be a recurrence sequence
of order $d$  of rational integers. Then, for every prime number $p$
and every positive integer $k$, the period of the sequence $(u_n)_{n
\ge 1}$ modulo $p^k$ is at most equal to $(p^d - 1) p^{k-1}$.

\pro
See Everest et al. \cite{EPSW03}, page 47. \cqfd

\proclaim Lemma 6.4.
Let $\alpha = [a_0; a_1, \ldots , a_{r-1}, b_0, b_1, \ldots , b_{s-1}, b_0, \ldots , b_{s-1}, \ldots]$
be a quadratic irrational number and denote by $(p_n / q_n)_{n \ge 0}$
the sequence of its convergents. Then, there exists an integer $t$ such that
$$
q_{n + 2s} - t q_{n+s} + (-1)^s q_n = 0
$$
for $n \ge r$. In particular, the sequence $(q_n)_{n \ge 0}$ satisfies
a linear recurrence with constant integral coefficients.

\pro
This result is included in the proof of Theorem 1 in \cite{LeSh93}. \cqfd

\vskip 5mm

\goodbreak

\centerline{\bf 6.2. Proofs}

\vskip 6mm

\noi {\it Preliminaries.}

Without any loss of generality, we consider real numbers in $(0,
1)$. We associate to every real irrational number $\alpha := [0;
a_1, a_2, \ldots]$ the infinite word ${\bf a} := a_1 a_2 \ldots $
formed by the sequence of the partial quotients of its fractional
part. Set
$$
p_{-1} = q_0 = 1, \quad p_0 = q_{-1} = 0,
$$
and
$$
{p_n \over q_n} = [0; a_1, \ldots, a_n], \quad
\hbox{for $n \ge 1$.}
$$
By the theory of
continued fractions, we have
$$
{q_n \over q_{n-1}} = [a_n; a_{n-1}, \ldots , a_1].
$$
This is one of the key tools of our proofs.

\bigskip

\noi {\it Proof of Theorems 2.3 and 3.3.}

Assume that the infinite word $a_{m+1} a_{m+2} \ldots $ is recurrent.
Then, there exists an increasing sequence of positive integers $(n_j)_{j \ge 1}$
such that
$$
\hbox{$a_{m+1} a_{m+2} \ldots a_{m+n_j}$
is a suffix of $a_{m+1} a_{m+2} \ldots a_{m+n_{j+1}}$, for $j \ge 1$.}
$$
Say differently, there are finite words $V_1, V_2, \ldots $ such that
$$
\hbox{$a_{m+1} a_{m+2} \ldots a_{m+n_{j+1}} =
V_j a_{m+1} a_{m+2} \ldots a_{m+n_{j}}$, for $j \ge 1$.}
$$
Actually, these properties are equivalent.

Let $\ell \ge 2$ be an integer.
Let $k \ge \ell^2 + 1$ be an integer.
By Dirichlet's {\it Schubfachprinzip}, there exist integers $i, j$
with $1 \le i < j \le k$ such that
$$
q_{m + n_i} \equiv q_{m + n_j} \pmod \ell , \quad
q_{m + n_i - 1} \equiv q_{m + n_j - 1} \pmod \ell
$$
and $j$ is minimal with this property.

Setting
$$
Q := |q_{m + n_i} q_{m+n_j-1}  - q_{m+n_i -1}   q_{m + n_j}|,
$$
we observe that
$$
\hbox{$\ell$ divides $Q$,}   \eqno (6.1)
$$
and we derive from Lemma 6.2 that
$$
\eqalign{
0 < Q & = q_{m + n_i} q_{m + n_j} \Bigl|
{q_{m + n_j - 1} \over q_{m + n_j}} - {q_{m + n_i -1} \over q_{m + n_i}} \Bigr|  \cr
& \le q_{m + n_i} q_{m + n_j} K(a_{m+n_i}, \ldots , a_{m+1})^{-2},  \cr}
$$
since the $n_i$ first partial
quotients of $q_{m + n_j - 1} / q_{m + n_j}$ and $q_{m + n_i -1}/q_{m + n_i}$
are the same, namely $a_{m+n_i}, \ldots , a_{m+1}$.  Furthermore, we have
$$
||Q \alpha|| \le || q_{m+n_i} (q_{m+n_j-1}\alpha)|| + ||
q_{m+n_i-1}(q_{m+n_j}\alpha)|| \le 2 q_{m + n_i} q_{m + n_j}^{-1}.
$$
Using that
$$
q_{m + n_i} \le 2 q_m K(a_{m+n_i}, \ldots , a_{m+1}),
$$
by Lemma 6.1, we finally get
$$
Q \cdot ||Q \alpha|| \le 8 q_m^2.   \eqno (6.2)
$$
It then follows from (6.1) and (6.2) that
$$
Q \cdot ||Q \alpha || \cdot |Q|_{\ell} \le 8 q_m^2 \ell^{-1},   \eqno (6.3)
$$
where $|Q|_{\ell}$ is equal to $\ell^{-a}$ if $\ell^a$ divides $Q$ but $\ell^{a+1}$
does not.
Since $\ell$ can be an arbitrary prime power, this proves Theorem 2.3.

Our proof shows that there are arbitrarily large integers $q$ such that
$$
q \ell \cdot ||q (\ell \alpha) || \le 8 q_m^2,
$$
which implies that
$$
c(\ell \alpha) \le { 8 q_m^2 \over \ell},
$$
and establishes Theorem 3.3. \cqfd

\bigskip

\noi {\it Proof of Corollary 2.4.}

Let ${\bf a}$ be an infinite Sturmian word.
We first claim that every prefix of finite length of ${\bf a}$ occurs
infinitely often in ${\bf a}$. Indeed, otherwise, there would exist a positive
integer $n$, a finite word $W$ and an infinite word ${\bf a}'$
such that ${\bf a} = W {\bf a}'$ and $p(n, {\bf a}') \le n$, which would imply
that ${\bf a'}$ is ultimately periodic, a contradiction with the assumption
that ${\bf a}$ is Sturmian.

Let  ${\bf a}$
be an infinite word on a finite alphabet ${\cal A}$ such
that there are positive integers $k$ and $n_0$ with
$$
p(n, {\bf a}) = n + k, \quad
\hbox{for $n \ge n_0$}.
$$
Then, by a result of  Cassaigne \cite{Cassa98},
there exist finite words $W, W_0, W_1$
on ${\cal A}$ and a Sturmian word
${\bf s}$ on $\{0, 1\}$ such that
$$
{\bf a} = W \phi({\bf s}),
$$
where $\phi({\bf s})$ denotes the infinite word obtained by
replacing in ${\bf s}$ every $0$ by $W_0$ and every $1$ by $W_1$. We
conclude by applying Theorem 2.3 with $m$ being the length of $W$.
\cqfd

\bigskip

\noi {\it Proof of Corollary 2.5.}

It is sufficient to note that, if $a_1 \ldots a_n$ and $a_1 \ldots
a_{n'}$ are palindromes with $n' > 2 n$, then $a_{n'-n+1} \ldots
a_{n'} = a_n \ldots a_1 = a_1 \ldots a_n$. The corollary then
follows from Theorem 2.3 applied with $m=0$. \cqfd

\bigskip

\noi {\it Proof of (1.3) when $\alpha$ is a quadratic irrationality
and ${\cal D}$ is bounded.}

Since ${\cal D}$ is bounded, every product $e_n = \prod_{1\le k\le n} d_k$
is divisible by a finite collection of prime numbers. Let $p_1,
\ldots , p_h$ be these primes and denote by $S$ the set of integers
which are divisible only by primes from $\{p_1, \ldots , p_h\}$. Let
$\alpha$ be a quadratic real number. By Lemma~6.4, the sequence
$(q_n)_{n \ge 0}$ of denominators of convergents to $\alpha$ is
eventually a recurrence sequence of positive integers. 
By Lemma 6.3, there exists a positive integer $C_1$
such that, for $i=1 , \ldots , h$ and $v \ge 1$,
the sequence $(q_n)_{n \ge 0}$ is eventually periodic modulo
$p_i^v$, with period length at most equal to $C_1 p_i^v$.

Consequently, there exists  a positive integer $C_2$ such that,
for every positive integer $\ell$
in $S$,  the sequence $(q_n)_{n \ge 1}$ modulo $\ell$
is eventually periodic of period at most $C_2 \ell$.

We need to slightly modify the proof of Theorem 2.3. Take $\ell =
e_n\in S$. Denote by $m$ the length of the preperiod of $(q_n)_{n
\ge 1}$ and by $d$ the length of the period of $(q_n)_{n \ge 1}$
modulo $\ell$. Observe that
$$
q_{m} \equiv q_{m +  d} \pmod \ell , \quad
q_{m + 1} \equiv q_{m + d + 1} \pmod \ell .
$$
We then set
$$
Q := |q_m q_{m+ d + 1} - q_{m + 1} q_{m + d}|
$$
and proceed exactly as in the proof of Theorem 2.3 to get that
$$
Q \cdot ||Q \alpha|| \le 2 q_{m+1}^2.
$$
Noticing that $|Q|_{\cal D}\le \ell^{-1}$ and $Q \le C_3^{\ell}$, for
some integer $C_3$ depending only on $p_1, \ldots , p_h$, this
establishes (1.3).  \cqfd

\bigskip

\goodbreak

\noi {\it Proof of Theorem 2.7.}

We keep the notation of the proof of Theorem 2.3.
By assumption, we can select a suitable sequence $(n_j)_{j \ge 1}$
with the property that $n_j < C_1^j$ for some integer $C_1 \ge 2$ and every $j \ge 1$.
Then, there are positive constants $C_2, C_3$, depending only
on $C_1$, such that
$$
Q \le C_2^{n_i + n_j},
$$
thus,
$$
\log \log Q \le C_3 \ell^2,
$$
since $i$ and $j$ are at most equal to $\ell^2+1$.
Combined with (6.3), this proves the theorem. \cqfd

\bigskip

\vskip 7mm

\vfill\eject

\centerline{\bf References}

\vskip 7mm

\beginthebibliography{999}

\bibitem{AlSh}
J.-P. Allouche and J. Shallit,
Automatic Sequences: Theory, Applications, Generalizations.
Cambridge University Press, Cambridge, 2003.

\bibitem{BaVe11}
D. Badziahin and S. Velani,
{\it Multiplicatively badly approximable numbers
and the mixed Littlewood conjecture},
Adv. Math. 228 (2011), 2766--2796.

\bibitem{BuLiv}
Y. Bugeaud,
Approximation by algebraic numbers,
Cambridge Tracts in Mathematics 160,
Cambridge, 2004.

\bibitem{Bu14}
Y. Bugeaud,
{\it Around the Littlewood conjecture in Diophantine approximation},
Publ. Math. Besan\c con Alg\`ebre Th\'eorie Nr. (2014), 5--18.

\bibitem{BDM}
Y. Bugeaud, M. Drmota, and B. de Mathan,
{\it On a mixed Littlewood conjecture in Diophantine approximation},
Acta Arith. 128 (2007), 107--124.

\bibitem{Cassa98}
J. Cassaigne,
{\it Sequences with grouped factors}.
In: DLT'97, Developments in
Language Theory III, Thessaloniki,
Aristotle University of Thessaloniki, 1998, pp. 211--222.

\bibitem{CaSw}
J. W. S. Cassels and H. P. F. Swinnerton-Dyer,
{\it On the product of three homogeneous linear forms and indefinite
ternary quadratic forms}, Philos. Trans. Roy. Soc. London,
Ser. A, 248 (1955), 73--96.

\bibitem{EFS}
M. Einsiedler, L. Fishman, and U. Shapira,
{\it Diophantine approximation on fractals},
Geom. Funct. Anal. 21 (2011), 14--35.

\bibitem{EKL}
M. Einsiedler, A. Katok, and E. Lindenstrauss,
{\it Invariant measures and the set of exceptions to the
Littlewood conjecture},
Ann. of Math. 164 (2006), 513--560.

\bibitem{EiKl07}
M. Einsiedler and D. Kleinbock,
{\it Measure rigidity and $p$-adic Littlewood-type problems},
Compositio Math. 143 (2007), 689--702.

\bibitem{EinLin}
M. Einsiedler and E. Lindenstrauss,
{\it On measures invariant under tori on quotients of semi-simple groups},
Ann. of Math. To appear.

\bibitem{EPSW03}
G. Everest, A. van der Poorten, I. Shparlinski, and T. Ward,
Recurrence sequences. Mathematical Surveys and Monographs, 104.
American Mathematical Society, Providence, RI, 2003.

\bibitem{LeSh93}
H. W. Lenstra and J. O. Shallit,
{\it Continued fractions and linear recurrences},
Math. Comp. 61 (1993), 351--354.

\bibitem{Lindenstrauss}
E. Lindenstrauss,
{\it Invariant measures and arithmetic quantum unique
ergodicity}, Ann. of Math. (2)  163  (2006),  165--219.

\bibitem{Lin10}
E. Lindenstrauss,
{\it Equidistribution in homogeneous spaces and number theory}.
In: Proceedings of the International Congress of Mathematicians. 
Volume I,  531--557, Hindustan Book Agency, New Delhi, 2010.

\bibitem{BdMTe}
B. de Mathan et O. Teuli\'e,
{\it Probl\`emes diophantiens simultan\'es},
Monatsh. Math. 143 (2004), 229--245.

\bibitem{MoHe38}
M. Morse and G. A. Hedlund,
{\it Symbolic dynamics},
Amer. J. Math. 60 (1938), 815--866.

\bibitem{MoHe40}
M. Morse and G. A. Hedlund,
{\it Symbolic dynamics II},
Amer. J. Math. 62 (1940), 1--42.

\bibitem{Per}
O. Perron,
Die Lehre von den Ketterbr\"uchen.
Teubner, Leipzig, 1929.

\bibitem{PoVe}
A. D. Pollington and S. Velani,
{\it On a problem in simultaneous Diophantine approximation:
Littlewood's conjecture},
Acta Math. 185 (2000), 287--306.

\vskip 4mm

\endthebibliography

\goodbreak

\vskip 6mm

\noindent Dmitry Badziahin  \hfill  Yann Bugeaud

\noindent   University of Durham  \hfill Universit\'e de Strasbourg

\noindent Department of Mathematical Sciences  \hfill Math\'ematiques

\noindent   South Rd  \hfill 7, rue Ren\'e Descartes

\noindent  Durham DH1 3LE  \hfill 67084 Strasbourg Cedex

\noindent   UK  \hfill France

\medskip

\noindent   {\tt dzmitry.badziahin@durham.ac.uk}  \hfill  {\tt bugeaud@math.unistra.fr}

\goodbreak

\bigskip

\noindent Manfred Einsiedler \hfill Dmitry Kleinbock

\noindent ETH Z\"urich, Departement Mathematik \hfill Brandeis University

\noindent R\"amistrasse 101 \hfill Department of Mathematics

\noindent CH-8092 Z\"urich \hfill Waltham, MA 02454

\noindent Switzerland \hfill USA

\medskip

\noindent {\tt manfred.einsiedler@math.ethz.ch} \hfill {\tt kleinboc@brandeis.edu}

\bye